%% file: mono_fin.tex
\documentclass[12pt]{article}
\usepackage[psamsfonts]{amsfonts}
\usepackage{amsmath,amssymb,amsthm}
\usepackage[dvips]{graphicx}
\usepackage{eepic}
\newtheorem{THM}{Theorem}[section]
\newtheorem{COR}[THM]{Corollary}

\newtheorem{LEM}[THM]{Lemma}
\newtheorem{PRP}[THM]{Proposition}

\input {def99c} 
\UseSection  
\usepackage[usenames]{color}

\newcommand{\ch}[1]{#1}
\newcommand{\ci}[1]{#1}
\newcommand{\ra}{\rightarrow}

\newcommand{\bigcupd} {\stackrel{\cdot}{\bigcup}}
\newcommand{\del}{\partial}

\newcommand{\shift}   {\!\!\!\!}

\newcounter{countC}  
\newcounter{countR}  
\newcommand{\prodstack}[2]{\prod_{\substack{#1 \\ #2}}}

\newcommand{\Z}{\Zbold}

\newcommand{\nn}{\nonumber}

\newcommand{\smallsup}[1] {{\scriptscriptstyle{({#1}})}}

\newcommand{\walk}{\vec{\omega}}
\newcommand{\ewalk}{\vec{\eta}}
\newcommand{\walkvec}[2]{\vec{\omega}^{\smallsup{#1}}_{#2}}
\newcommand{\walkcoor}[2]{\omega^{\smallsup{#1}}_{#2}}

\newcommand{\mc}[1]{\mathcal{#1}}
\newcommand{\mP}{\mathbb{P}}
\newcommand{\mQ}{\mathbb{Q}}
\newcommand{\mE}{\mathbb{E}}
\newcommand{\cp}{\overset{\mP}{\ra}}
\newcommand{\cas}{\overset{a.s.}{\ra}}
\newcommand{\eqn}[1]{\eq #1 \en}
\newcommand{\sN}{{\sss N}}
\newcommand{\dimres}{\ch{9}}
\newcommand{\dimlow}{\ci{8}}

\title  {
        Monotonicity for excited random walk\\ in high dimensions
        }

\author{
Remco van der Hofstad\footnote{Department of Mathematics and
Computer Science, Eindhoven University of Technology, P.O.\ Box
513, 5600 MB Eindhoven, The Netherlands. E-mail {\tt
rhofstad@win.tue.nl}}
\and
Mark Holmes\footnote{Department of Statistics,
The University of Auckland, Private Bag 92019, Auckland 1142,
New Zealand. E-mail {\tt mholmes@stat.auckland.ac.nz}}
}

\begin{document}

\maketitle

    \begin{abstract}
    We prove that the drift $\theta(d,\beta)$ for excited random walk
    in dimension $d$ is monotone in the excitement parameter $\beta\in [0,1]$, when $d\ge \dimres$.
    \end{abstract}


\section{Introduction}
\label{sec-intro}
In this paper we study {\it excited random walk},  where the random walker has a drift
in the direction of the first component each time the walker visits a
new site. It was shown that this process has ballistic behaviour when $d\geq 2$ in
\cite{BW03,Kozm03, Kozm05}, while there is no ballistic behaviour in one dimension.  \ci{A nontrivial strong law of large numbers can then be obtained for $d\ge2$ using renewal techniques (see for example \cite{SZ99}, \cite{Zern05}).
In \cite{HH07} a perturbative expansion was introduced and used to prove a weak law of
large numbers and a central limit theorem for excited random walk in
dimensions $d>5$ and $d>8$ respectively, with sufficiently small excitement
parameter.  More recently \cite{BR07} explicitly proved a SLLN and established a functional central limit theorem in dimensions $d\ge 2$.}

Included in \cite{HH07} is an explicit representation of the drift in
terms of the expansion coefficients.  In this paper we use this
representation, together with improved diagrammatic estimates to prove
that in dimensions $d\ge \dimres$, the drift for
excited random walk is (strictly) increasing in the excitement parameter $\beta$.

\subsection{Main results}
\label{sec-results}
The main result of this paper is the following theorem.

\begin{THM}[\ch{Monotonicity of the speed}]
\label{thm:main}
For all $d\ge \dimres$, and $\beta \in [0,1]$, the drift for excited random walk
in dimension $d$ with excitement parameter $\beta$ is strictly increasing in $\beta$.
\end{THM}

We are also able to show that for $d\ge \dimlow$, there exists
$\beta_0(d)$ such that the drift for ERW is strictly increasing in $\beta\in [0,\beta_0]$.

Simulations \cite{BR07} suggest that the limiting variance of the first coordinate
is {\it not} monotone in the excitement parameter $\beta$ in 2 dimensions.
We expect that using the approach introduced in this paper we
can show that the variance is monotone decreasing in $\beta$
when the dimension is taken sufficiently high.  By \cite{HH07}, the
variance of the first coordinate is equal to $\sigma_{\beta}^2 n(1+o(1))$
for some asymptotic variance $\sigma_{\beta}^2$,
and based on our methods, we expect that
$\sigma_{\beta}^2=d^{-1}-d^{-2}\beta^2+\beta^2\mc{O}\left(d^{-3}\right)$,
showing that, in sufficiently high $d$, $\beta\mapsto \sigma_{\beta}^2$ is decreasing.

Although we only consider the {\em once}-excited random walk in
this paper, the general {\em multi}-excited random walk can be
handled with very minor modifications, yielding a result at
least as strong as Theorem \ref{thm:main}.
A large part of the methodology in this paper can be applied
more generally.  Given the present context, a natural example
is a random walk in an environment that is random \ci{in the first few
coordinates only}, with the expected drift induced by the environment
denoted by $\vec{\beta}$. Laws of large numbers and functional
central limit theorems have been proved for general versions of such
random walks in random environment in \cite{BolSznZei03}.  
\ci{We intend to investigate monotonicity 
issues for the speed in such models in the near future.}

We first introduce some notation. A nearest-neighbour
random walk path $\vec{\eta}$ is a sequence
$\{\eta_i\}_{i=0}^\infty$ for which $\eta_i\in \Z^d$
and $\eta_{i+1}-\eta_i$ is a nearest-neighbour of the origin
for all $i\geq 0$.
For a general nearest-neighbour path $\vec{\eta}$ with $\eta_0=0$, we write
$p^{\vec{\eta}_i}(x_i, x_{i+1})$
for the conditional probability that the walk steps from $\eta_i=x_i$ to
$x_{i+1}$, given the history of the path $\vec{\eta}_i=(\eta_0, \ldots, \eta_i)$.
We write $\walk_n$ for the $n$-step path of excited random walk (ERW),
and $\Qbold$ for the law of $\{\walk_n\}_{n=0}^{\infty}$, i.e.,
for every $n$-step nearest-neighbour path $\ewalk_n$,
    \eq
    \lbeq{SIRP}
    \Qbold(\walk_n=\vec{\eta}_n)
    =\prod_{i=0}^{n-1} p^{\vec{\eta}_i}(\eta_{i},\eta_{i+1}),
    \en
where, for $i=0$, $p^{\varnothing}(0,\eta_{1})$ is the probability to jump to
$\eta_1$ in the first step, and
    \eq
    p^{\vec{\eta}_i}(\eta_{i},\eta_{i+1})=p_0(\eta_{i+1}-\eta_i)\delta_{\eta_i, \vec{\eta}_{i-1}}
    +p_\beta (\eta_{i+1}-\eta_i)[1-\delta_{\eta_i, \vec{\eta}_{i-1}}],
    \en
where $\delta_{\eta_i, \vec{\eta}_{i-1}}$ denotes the indicator that $\eta_i=\eta_j$ for some
$0\leq j\leq i-1$, and where, for $\beta\in [0,1]$,
    \eq
    p_\beta(x) =\frac{1+\beta e_1\cdot x}{2d} I[|x|=1].
    \en
Here $e_1=(1,0, \ldots, 0)$, and $x\cdot y$ is the
inner-product between $x$ and $y$.  We will usually drop
the indicator function here, and leave it implicit in the
notation that our walks take nearest-neighbour steps.  In words, the random walker
gets excited each time he/she visits a new site, and
when the random walk is excited, it has a positive drift
in the direction of the first coordinate. For a description
in terms of cookies, see \cite{Zern05}. We let $\Ebold$
denote expectation with respect to $\Qbold$.

It is known that in dimensions $d \ge 2$, excited random walk
has a positive drift $\theta=\theta(\beta,d)$ satisfying
$n^{-1}\omega_n\cas \theta$ and that a (functional) central
limit theorem holds \cite{BW03, BR07, HH07, Kozm03, Kozm05}.
For $d=1$, it is known that ERW is recurrent and diffusive
\cite{Dav99} except in the trivial case $\beta=1$.  Additional
results on one-dimensional (multi)-excited random walks can be
found in \cite{AR05, Zern05}.

\section{\ch{An overview of the proof and the expansion}}
\label{sec-fromexpansion}
In this section we recall some results and notation from \cite{HH07}.
If $\ewalk$ and $\walk$ are two paths of length at least $j$ and $m$
respectively and such that $\eta_j=\omega_0$, then the concatenation
$\ewalk_j\circ \walk_m$ is defined by
    \eq
    \lbeq{concat}
    (\ewalk_j\circ \walk_m)_{\ch{i}}=\left\{
    \begin{array}{lll}
    &\eta_i&\text{when }0\le i\leq j,\\
    &\omega_{i-j} &\text{when }j \leq i \leq m+j.
    \end{array}\right.
    \en
Given $\ewalk_m$, we define a probability measure $\Qbold^{\ewalk_m}$
on walks path starting from $\eta_m$ by specifying its value on particular
cylinder sets (in a consistent manner) as follows
    \eq
    \lbeq{cylinders}
    \Qbold^{\ewalk_m} (\walk_n=\vec{\mu}_n)
    \equiv\prod_{i=0}^{n-1} p^{\ewalk_m\circ \ci{\vec{\mu}_i}}(\mu_{i},\mu_{i+1}),
    \en
and extending the measure to all finite-dimensional cylinder sets in the natural (consistent) way.  \ci{When $\Qbold(\walk_m=\ewalk_m)>0$, \refeq{cylinders} is also $\Qbold(\walk_{m+n}=\ewalk_m\circ \vec{\mu}_n|\walk_m=\ewalk_m)$.}  We write $\Ebold^{\ewalk_m}$ for the expected value with respect to $\Qbold^{\ewalk_m}$. \ci{In \cite{HH07}, a perturbative expansion was derived for the two-point function $c_n(x)=\Qbold(\omega_n=x)$, giving rise to a recursion relation of the form
    \eqalign
    c_{n+1}(x)=\sum_y p^{\varnothing}(0,y)c_n(x-y)+\sum_{m=2}^{n+1}\sum_y \pi_m(y)c_{n+1-m}(x-y).\lbeq{recrel}
    \enalign
This expansion was used to prove a law of large numbers
and central limit theorem for ERW.  We next discuss the coefficients $\pi_m(y)$ and some results of this expansion.}

\paragraph{The expansion coefficients.} \ci{The lace expansion coefficients involve the following factors}.  For $N\geq 1$, let
    \eq
    \lbeq{DeltaNdef}
    \Delta^{\smallsup{N}}_{j_{N}+1}= \big(p^{\walkvec{N-1}{j_{N-1}+1}\circ
    \walkvec{N}{j_{N}}}-p^{\walkvec{N}{j_{N}}}\big)(\walkcoor{N}{j_{N}},\walkcoor{N}{j_{N}+1}),
    \en
with $j_0\equiv 0$.  The difference \refeq{DeltaNdef} is identically zero when the histories
$\walkvec{N-1}{j_{N-1}+1}\circ \walkvec{N}{j_{\sN}}$ and $\walkvec{N}{j_{\sN}}$ give the same transition
probabilities to go from $\walkcoor{N}{j_{\sN}}$ to $\walkcoor{N}{j_{\sN}+1}$.  For excited random walk, $\Delta^{\smallsup{N}}_{j_{\sN}+1}$ is non-zero precisely when $\walkcoor{N}{j_{\sN}}$ has already been visited by $\walkvec{N-1}{j_{N-1}+1}$ but not by $\walkvec{N}{j_{N}-1}$, so that
    \eqalign
    \lbeq{Deltabound}
    |\Delta^{\smallsup{N}}_{j_{N}+1}|=&\left|\frac{\beta e_1\cdot(\walkcoor{N}{j_{N}+1}-\walkcoor{N}{j_{N}})}{2d}\left[I[\walkcoor{N}{j_{N}}\notin \walkvec{N-1}{j_{N-1}}\circ \walkvec{N}{j_{\sN}-1}]-I[\walkcoor{N}{j_{N}}\notin \walkvec{N}{j_{\sN}-1}]\right]\right|\\
    \le &\frac{\beta}{2d}I[\walkcoor{N}{j_{N}+1}=\walkcoor{N}{j_{N}}\pm e_1]I[\walkcoor{N}{j_{N}}\in \walkvec{N-1}{j_{N-1}}\setminus \walkvec{N}{j_{\sN}-1}]\le \frac{\beta}{2d}I[\walkcoor{N}{j_{N}+1}=\walkcoor{N}{j_{N}}\pm e_1]I[\walkcoor{N}{j_{N}}\in \walkvec{N-1}{j_{N-1}}].\nn
    \enalign

Define $\mc{A}_{m,\sN}=\{(j_1, \dots, j_{\sN})\in \Z_+^{\sN}:\sum_{l=1}^Nj_l=m-N-1\}$, $\mc{A}_{\sN}=\bigcupd_{m}\mc{A}_{m,\sN}$ and
    \eqalign
    &\pi_m^{\smallsup{N}}(x,y)\lbeq{piNxydef}\\
    =&\shift
    \sum_{\vec{j}\in \mc{A}_{m,N}}\Ebold_{\sss 0}^{\varnothing}\Bigg[\sum_{\walkvec{0}{1}}
    \Ebold_{\sss 1}^{\walkvec{0}{1}}\Big[\sum_{\walkcoor{1}{j_1+1}}\Delta^{\smallsup{1}}_{j_1+1}~
\Ebold_{\sss 2}^{\walkvec{1}{j_1+1}}
   \big[\sum_{\walkcoor{2}{j_2+1}}\Delta^{\smallsup{2}}_{j_2+1}~
   \cdots \Ebold_{\sss N}^{\walkvec{N-1}{j_{N-1}+1}}
   \big[\sum_{\walkcoor{N}{j_{\sN}+1}}\Delta^{\smallsup{N}}_{j_{N}+1}I_{\{\omega^{(N)}_{j_{\sN}}=x, \omega^{(N)}_{j_{\sN}+1}=y\}}\big] \cdots \Big]\Big]\Bigg].
   \nn
   \enalign
Then we define
    \eq
    \pi_m(x,y)=\sum_{N=1}^{\infty} \pi_m^{\smallsup{N}}(x,y),  \quad
    \pi^{\smallsup{N}}(x,y)=\sum_{m}\pi_m^{\smallsup{N}}(x,y),
    \quad  \text{and} \quad \pi_m(y)=\sum_{N=1}^{\infty}\sum_x \pi_m^{\smallsup{N}}(x,y).
    \lbeq{piotherdef}
    \en
Note that the quantities $\pi_m^{\smallsup{N}}$ are all zero when $N+1>m$,
and that all of the above quantities depend on $\beta$.  We emphasize that,
conditionally on $\walkvec{M}{j_{M}+1}$,
the probability measure $\ch{\Qbold}_{\sss M+1}^{\walkvec{M}{j_{M}+1}}$ is the law of
$\walkvec{M+1}{j_{M+1}+1}$, \ch{i.e.,} $\walkvec{M}{j_{M}+1}$
acts as the history for $\walkvec{M+1}{j_{M+1}+1}$.

In \cite{HH07}, it is also
shown that if $\lim_{n \ra \infty}\sum_{m=2}^n \sum_x x\pi_m(x)$ exists
and $n^{-1}\omega_n\cp \theta$, then
    \eqalign
    \theta(\beta,d)& = \sum_x xp^{\varnothing}(0,x) +\sum_{m=2}^\infty \sum_x x\pi_m(x)\lbeq{theta}.
    \enalign
\ci{\paragraph{Strategy of the proof of Theorem \ref{thm:main}.}  We shall explicitly differentiate
the right hand side of \refeq{theta}, and prove that this derivative
is positive for all $\beta\in [0,1]$, when }$d\ge \dimres$.
From \refeq{theta} and using the fact that $\sum_y \pi_m(x,y)=0$
(recall \refeq{piotherdef}),
we have
    \eqalign
    \sum_y  y\pi_m(y)=&\sum_{x,y}  (y-x)\pi_m(x,y),\lbeq{xy}
    \enalign
so that
    \eqalign
    \lbeq{theta2}
    \theta(\beta, d)=&\frac{\beta e_1}{d}+\sum_{m=2}^\infty\sum_{N=1}^{\infty}\sum_{x,y}  (y-x)\pi_m^{\smallsup{N}}(x,y).
    \enalign
Letting $\varphi_m^{\smallsup{N}}(x,y)=\frac{\del}{\del\beta}\pi_m^{\smallsup{N}}(x,y)$
and assuming that the limit can be taken through the infinite sums, we then have
    \eqalign
    \lbeq{deriv1}
    \frac{\del \theta}{\del \beta}(\beta,d)&=\frac{e_1}{d}+\sum_{N=1}^{\infty}\sum_{m=2}^\infty\sum_{x,y}  (y-x)\varphi_m^{\smallsup{N}}(x,y).
    \enalign
Since $\varphi_m^{\smallsup{N}}(x,y)\equiv 0$ unless $|x-y|=1$, we have that
    \eqalign
    \left|\frac{\del \theta}{\del \beta}(\beta,d)-\frac{e_1}{d}\right|\le \sum_{N=1}^{\infty}\sum_{m=2}^\infty\sum_{x,y}|\varphi_m^{\smallsup{N}}(x,y)|.\lbeq{needed1}
    \enalign
We conclude that $\frac{\del \theta_{1}}{\del \beta}(\beta,d)$, which is the first coordinate of $\frac{\del \theta}{\del \beta}(\beta,d)$,
is positive for any $\beta$ at which $\sum_{N=1}^{\infty}\sum_{m=2}^\infty\sum_{x,y}|\varphi_m^{\smallsup{N}}(x,y)|<d^{-1}$.
This is what we shall prove in the remainder of this paper, which is organised
as follows. In Section \ref{sec-pibd}, we start by proving bounds on
$\pi_m^{\smallsup{N}}$. These bounds will be crucially used to
prove bounds on $\varphi_m^{\smallsup{N}}$ in Section \ref{sec-derivpi}.
The results in Section \ref{sec-derivpi} are used in Section \ref{sec-pfmainthm}
to prove Theorem \ref{thm:main}.


\section{Bound on $\pi$}
\label{sec-pibd}
Before proceeding to the proof of Theorem \ref{thm:main}, we prove
a new bound on $\sum_{x,y}\sum_m|\pi^{\smallsup{N}}_m(x,y)|$.
The proof of this new bound makes use of Lemmas \ref{lem:decomposition}
and \ref{lem:togreens} below.  For the first of these lemmas we need
to introduce some notation.

Let $f_{i,{j}_i}(\walkvec{i-1}{m},\walkvec{i}{j_i})\ge0$, $i=0,\dots,N,$ be
measurable functions from the set of (ordered) pairs of finite random walk
paths $(\walkvec{i-1}{m},\walkvec{i}{j_i})$ such that $m<\infty$ and
$\walkcoor{i-1}{m}=\walkcoor{i}{0}$ (the former is defined to be the
origin if $\walkvec{i-1}{m}=\varnothing$).
Recall that $\mE^{\walkvec{i-1}{m}}$ denotes expectation with respect to the law of a self-interacting random walk (ERW in this paper)
$\walkvec{i}{}$ with given (finite)
history $\walkvec{i-1}{m}$ (i.e., conditional on the first $m$
steps of the walk being $\walkvec{i-1}{m}$).  We write
$\mE^{l,\walkvec{i-1}{m}}$ to distinguish expectation with respect to different laws
(indexed by $l$), i.e., if $l_1\ne l_2$ then $\mQ^{l_1,\walkvec{i-1}{m}}$
and $\mQ^{l_2,\walkvec{i-1}{m}}$ may be different self-interacting random walk
laws (with the same given history).

Given $\vec{f}_{\sN}=(f_{1,j_1},\dots,f_{N,j_{\sN}})$ and $k=0, \ldots, N$, we define
    \eqalign
    \Pi_{\sN}^{\sss(k)}(\vec{f}_{\sN})\equiv & \sum_{\vec{j}\in \mc{A}_{\sN}}
    \mE^{\varnothing}_0\Bigg[f_{0,j_0}(\varnothing,\walkvec{0}{j_0})
    \mE^{\walkvec{0}{j_0+1}}_1\big[f_{1,{j}_1}(\walkvec{0}{j_0+1},\walkvec{1}{j_1})\dots \sum_{l=1}^{j_k}\mE^{l,\walkvec{k-1}{j_{k-1}+1}}_{k}[f_{k,{j}_k}(\walkvec{k-1}{j_{k-1}+1},\walkvec{k}{j_k})
    \dots\nn\\ &\qquad \qquad \mE^{\walkvec{N-1}{j_{N-1}+1}}_{N}
    [f_{N,{j}_{\sN}}(\walkvec{N-1}{j_{N-1}+1},\walkvec{N}{j_{\sN}})]\dots\big]\Bigg].
    \lbeq{bigaltpidef}
    \enalign
We further let $\Pi_{\sN}(\vec{f}_{\sN})$ be identical to $\Pi_{\sN}^{\sss(k)}(\vec{f}_{\sN})$,
apart from the fact that $\sum_{l=1}^{j_k}\mE^{l,\walkvec{k-1}{j_{k-1}+1}}_{k}$
is replaced with $\mE^{\walkvec{k-1}{j_{k-1}+1}}_{k}$. A crucial ingredient in
obtaining bounds on lace expansion coefficients is the following result:

\begin{LEM}[\ci{Recursive bounds for expansion coefficients}]
\label{lem:decomposition}
Let $\walkvec{0}{}, \dots, \walkvec{N}{}$ be any collection of
$N$ self-interacting random walks defined on the same
probability space $(\Omega,\mc{F},\mP)$.  Suppose that
$f_{i,{j}_i}\ge0$, $i=0,\dots,N$ are such that for each $i=0, \dots, N$
there exist constants $K_i\ge 0$, and functions $\kappa_i\ge 0$
(with $\kappa_{-1}\equiv 1$) such that
    \eqalign
    \sum_{j_i=0}^{\infty}\kappa_i(j_{i}) \mE^{\walkvec{i-1}{m}}
    [f_{i,{j}_i}(\walkvec{i-1}{m},\walkvec{i}{j_i})]\le K_i\kappa_{i-1}(m),\lbeq{piecebound}
    \enalign
for each $m$, uniformly in $\walkvec{i-1}{m}$.  Then
    \eqalign
    \lbeq{PiNbd}
    \Pi_{\sN}(\vec{f}_{\sN}) \le \prod_{i=0}^NK_i.
    \enalign
The conclusion in \refeq{PiNbd} also holds for $\Pi_{\sN}^{\sss(k)}(\vec{f}_{\sN})$
if there exist $K_i, \kappa_i\ge 0$ such that \refeq{piecebound}
holds for $i \ne k$, and for $i=k$,
    \eqalign
    \sum_{j_k=0}^{\infty}\kappa_k(j_{k})\sum_{l=1}^{j_k} \mE^{l,\walkvec{k-1}{m}}[f_{k,{j}_k}(\walkvec{k-1}{m},\walkvec{k}{j_k})]\le K_i\kappa_{k-1}(m),\lbeq{altpiecebound}
    \enalign
for each $m$, uniformly in $\walkvec{i-1}{m}$.
\end{LEM}

\proof Taking the sum over $j_{\sN}$ inside the first $N-1$ expectations
and then using \refeq{piecebound} with $i=N$ gives a factor
$K_{\sN}\kappa_{N-1}(j_{N-1})$.  The $K_{\sN}$ can be taken outside all
of the expectations and sums, while the $\kappa_{N-1}(j_{N-1})$
remains inside the sum over $j_{N-1}$.  Proceeding inductively
using \refeq{piecebound} we obtain the first result.  The proof
of the second result is identical except that at some point we
use \refeq{altpiecebound} instead of \refeq{piecebound}.
\qed

\vspace{.5cm}

Let $\mP_{d}$ denote the law of simple symmetric random walk in
$d$ dimensions, beginning at the origin, and let $D_d(x)=I[|x|=1]/(2d)$ be the
simple random walk step distribution. We will make use of the
convolution of functions, which is
defined for absolutely summable functions $f,g$ on $\Zd$ by
    \eq
    (f*g)(x) = \sum_y f(y) g(x-y).
    \en
Let $f^{*k}(x)$ denote the
$k$-fold convolution of $f$ with itself, and let
$G_d(x)=\sum_{k=0}^{\infty} D_d^{*k}(x)$ denote the Green's
function for this random walk. We shall sometimes make use of the
representation
    \eqn{
    \lbeq{convGd}
    G_d^{*i}(x)=\sum_{k=0}^{\infty}\sum_{\vec{m}_i:m_1+\dots+m_i=k} D_d^{*(m_1+\dots +m_i)}(x)=\sum_{k=0}^{\infty}\frac{(k+i-1)!}{(i-1)!k!} \mP_d(\omega_k=x),\quad \text{for }i\geq 1.
    }

For $i\geq 0$, let
    \eqalign
    \mc{E}_i(d)=&
    \sup_{v\in \Z^{d-1}}
    \Big(\big(\frac{d}{d-1}\big)^{i+1}G_{d-1}^{*(i+1)}(v)-\delta_{0,v}\Big).\lbeq{Eidef}
    \enalign

\begin{LEM}[\ch{Diagrammatic bounds for ERW}]
\label{lem:togreens}
For excited random walk, uniformly in $u\in \Z^d$, and, for $i\geq 0$,
    \eqalign
    \sum_{j=0}^{\infty}\frac{(j+i)!}{j!}\mQ^{\vec{\eta}_m}(\omega_{j}=u)\le &i!\big(\frac{d}{d-1}\big)^{i+1}G_{d-1}^{*(i+1)}(0),\lbeq{Gbound1}\\
    \sum_{j=1}^{\infty}\frac{(j+i)!}{j!}\mQ^{\vec{\eta}_m}(\omega_{j}=u)
    \le &i!\mc{E}_i(d)\lbeq{Gbound2}.
    \enalign
\end{LEM}

\proof
Let $j-\mc{N}_j$ be the number of steps that the walk $\vec{\omega}_j$ takes in the first coordinate.  Observe that independently of $\vec{\eta}$, $\mc{N}_j\sim Bin(j,q_d)$, where $q_d=(d-1)/d$.
If we consider $\vec{\omega}_j$ as the initial position and first $j$ steps of an infinite walk $\vec{\omega}$, then the sequence $\{\mc{N}_j\}_{j\ge 0}$ is a random walk on $\mathbb{Z}_+$ taking i.i.d.\ steps that are either +1 or 0 with probability $\ci{q_d}$ and $1-q_d$ respectively.  The random time that such a walk spends at any level $l$ has a Geometric distribution with parameter $q_d$.
Thus we obtain that, for every $i\geq 0$, and writing $\mc{P}$ for the law of $\{\mc{N}_j\}_{j=0}^{\infty}$,
    \eqalign
    \frac{(j+i)!}{j!}\mc{P}(\mc{N}_j=l)
    &= \frac{(j+i)!}{j!}\frac{j!}{l!(j-l)!} q_d^l (1-q_d)^{j-l}
    =q_d^{-i} \frac{(l+i)!}{l!}\mc{P}(\mc{N}_{j+i}=l+i),\nn
    \enalign
so that, for $m\le l$,
    \eqalign
    \sum_{j=m}^{\infty}  \frac{(j+i)!}{j!}\mc{P}(\mc{N}_j=l)=
    q_d^{-(i+1)}  \frac{(l+i)!}{l!}.\lbeq{tailwalk}
    \enalign

To prove \refeq{Gbound1} note that
    \eqalign
    \sum_{j=0}^{\infty}\frac{(j+i)!}{j!}\mQ^{\vec{\eta}_m}(\omega_{j}=u)=
    &\sum_{j=0}^{\infty}\frac{(j+i)!}{j!}\sum_{l=0}^j\mQ^{\vec{\eta}_m}(\omega_{j}=u|\mc{N}_j=l)
    \mc{P}(\mc{N}_j=l)\nn\\
    \le&\sum_{l=0}^\infty\mP_{d-1}(\omega_{l}=u^{[2,\dots,d]}-\eta_{m}^{[2,\dots,d]})
    \sum_{j=l}^{\infty}\frac{(j+i)!}{j!}\mc{P}(\mc{N}_j=l)\nn\\
    \le& q_d^{-(i+1)}\sup_{v\in \Z^{d-1}}\sum_{l=0}^\infty
    \mP_{d-1}(\omega_{l}=v)\frac{(l+i)!}{l!}\lbeq{cool1}.
    \enalign

By \refeq{convGd},
\refeq{cool1} is equal to $i! q_d^{-(i+1)}\sup_{v\in \Z^{d-1}}G_{d-1}^{*(i+1)}(v)$.
By \cite[Lemma B.3]{HS92b}, the supremum occurs at $v=0$. Using also that $q_d^{-1}=d/(d-1)$,
this proves \refeq{Gbound1}.

The bound \refeq{Gbound2} is proved similarly.  For $i\geq 0$ we can write
    \eqalign
    \sum_{j=1}^{\infty}\frac{(j+i)!}{j!}\mQ^{\vec{\eta}_{m}}
    (\omega_{j}=u)\le&\sup_{v\in \Z^{d-1}}\sum_{l=0}^\infty
    \mP_{d-1}(\omega_{l}=v)\sum_{j=l \vee 1}^{\infty}\frac{(j+i)!}{j!}
    \mc{P}(\mc{N}_j=l)\nn\\
    =&\sup_{v\in \Z^{d-1}}\Big(\sum_{l=0}^\infty\mP_{d-1}(\omega_{l}=v)
    \Big[\sum_{j=l}^{\infty}
    \frac{(j+i)!}{j!}\mc{P}(\mc{N}_j=l)
    -\delta_{0,l}i!\mc{P}(\mc{N}_0=0)\Big]\Big)\nn\\
    =&\sup_{v\in \Z^{d-1}}\Big(q_d^{-(i+1)}\sum_{l=0}^\infty
    \frac{(l+i)!}{l!}\mP_{d-1}(\omega_{l}=v)
    -i!\delta_{0,v}\Big)\nn\\
    =&i!\sup_{v\in \Z^{d-1}}\big(q_d^{-(i+1)}G_{d-1}^{*(i+1)}(v)-\delta_{0,v}\big),
    \lbeq{cool2}
    \enalign
since $\mc{P}(\mc{N}_0=0)=1$ and
$\sum_{l=0}^\infty \mP_{d-1}(\omega_{l}=v) \delta_{0,l}=\delta_{0,v}$,
and following the steps in \refeq{cool1} above.
\qed

Define
    \eqalign
    a_d=\frac{d}{(d-1)^{2}}G^{*2}_{d-1}(0)\lbeq{addef}.
    \enalign
\begin{PRP}[Bounds on the expansion coefficients]
\label{prp:pibound}
For all $N\ge 1$,
\eqalign
\sum_{x,y}\sum_m|\pi_m^{\smallsup{N}}(x,y)|\le
\begin{cases}
\beta d^{-1}\mc{E}_0(d) & N=1,\\
\beta^Nd^{-1}(d-1)^{-1}G_{d-1}(0)\mc{E}_1(d)\left(a_d\right)^{(N-2)I_{\{N>2\}}} & N>1.\lbeq{pibound}
\end{cases}
\enalign
\end{PRP}

Let
    \eqalign
    f_{0,j_0}(\vec{\eta}_m,\vec{\omega}_{j_0})=&I_{\{j_0=0\}},\quad f_{1,j_1}(\vec{\eta}_m,\vec{\omega}_{j_1})=\frac{\beta}{2d}I_{\{\omega_{j_1}= \eta_{0}\}}\sum_{\walkcoor{1}{j_1+1}}I_{\{\omega_{j_1+1}= \omega_{j_1}\pm e_1\}}I_{\{j_1 \text{ is odd}\}},\nn\\
    f_{i,j_i}(\vec{\eta}_m,\vec{\omega}_{j_i})=&\frac{\beta}{2d}I_{\{\omega_{j_i}\in \vec{\eta}_m\setminus \vec{\omega}_{j_i-1}\}}\sum_{\omega_{j_i+1}}I_{\{\omega_{j_i+1}= \omega_{j_i}\pm e_1\}}, \quad \text{ for } i>1.\lbeq{fdefs}
    \enalign
We will use Lemma \ref{lem:decomposition} for excited random walks, together with the following lemma to prove Proposition \ref{prp:pibound}.
    \begin{LEM}[Ingredients for bounding the coefficients for ERW]
    \label{lem:fbounds}
    For excited random walks with $f_{i,j_i}$ defined in \refeq{fdefs},
    \eqalign
    &\sum_{j_0=0}^{\infty}(j_0+1)\mE^{\vec{\eta}_m}_0[f_{0,{j}_0}]\le 1,\nn\\ &\sum_{j_1=0}^{\infty}(j_1+1)\mE^{\vec{\eta}_m}_1[f_{1,{j}_1}]\le \frac{\beta}{d}\mc{E}_1(d), \qquad \sum_{j_1=0}^{\infty}\mE^{\vec{\eta}_m}_1[f_{1,{j}_1}]\le \frac{\beta}{d}\mc{E}_0(d),\nn\\
    &\sum_{j_i=0}^{\infty}(j_i+1)\mE^{\vec{\eta}_m}_i[f_{i,{j}_i}]\le m\beta a_d, \quad i=2, \dots, N-1\nn\\
    &\sum_{j_{\sN}=0}^{\infty}\mE^{\vec{\eta}_m}_{\sN}[f_{\sN,{j}_{\sN}}]\le m\frac{\beta}{d-1}G_{d-1}(0).\lbeq{fbounds}
    \enalign
    \end{LEM}
\proof
The first bound is trivial.
For the second bound, since the conditions that $j\ge 0$
and $j$ is odd imply that $j\ge 1$, we have
    \eqalign
    \sum_{j=0}^{\infty}(j+1)\mE^{\vec{\eta}_m}_1[f_{1,{j}}]\le & \frac{\beta}{2d}\sum_{j=1}^{\infty}(j+1)\mE^{\vec{\eta}_m}[I_{\{\omega_{j}= \eta_{0}\}}\sum_{\walkcoor{1}{j+1}}I_{\{\omega_{j+1}= \omega_{j}\pm e_1\}}]\nn\\
    =&\frac{\beta}{d}\sum_{j=1}^{\infty}(j+1)\mQ^{\vec{\eta}_m}(\omega_{j}= \eta_{0})\le \frac{\beta}{d}\mc{E}_1(d),
    \enalign
where the last inequality holds by \refeq{Gbound2} with $i=1$.  Similarly
\refeq{Gbound2} with $i=0$ gives us the third bound.

For the fourth bound, using
    \eqalign
    I_{\{\omega_{j_i}\in \vec{\eta}_m\setminus \vec{\omega}_{j_i-1}\}}\le \sum_{l=0}^{m-1}I_{\{\omega_{j_i}=\eta_l\}}\nn
    \enalign
and proceeding as for the second bound we see that
    \eqalign
    \sum_{l=0}^{m-1}\sum_{j_i=0}^{\infty}(j_i+1)\mE^{\vec{\eta}_m}_i[f_{i,{j}_i}]\le &\sum_{l=0}^{m-1}\frac{\beta}{d}\sum_{j=0}^{\infty}(j+1)\mQ^{\vec{\eta}_m}(\omega_{j}= \eta_{l})\le m\beta a_d,\nn
    \enalign
where we have used \refeq{Gbound1} with $i=1$ in the last step.

For the last bound, note that
    \eqalign
    \sum_{j=0}^{\infty}\mE^{\vec{\eta}_m}_{\sN}[f_{\sN,j}]
    =&\frac{\beta}{2d}\sum_{j=0}^{\infty}\mE^{\vec{\eta}_m}[I_{\{\omega_{j}\in \vec{\eta}_{\ch{m}}\setminus \vec{\omega}_{j-1}\}}\sum_{\omega_{j+1}}I_{\{\omega_{j+1}= \omega_{j}\pm e_1\}}]=\frac{\beta}{d}\sum_{j=0}^{\infty}\mQ^{\vec{\eta}}(\omega_{j}\in \vec{\eta}_m\setminus \vec{\omega}_{j-1})\lbeq{flastb}\\
    \le&\sum_{k=0}^{m-1}\frac{\beta}{d}\sum_{j_{\sN}=0}^{\infty}\mQ^{\vec{\eta}_m}(\omega_{j_{\sN}}=\eta_k)
    \le \frac{\beta}{d}m\sup_u\sum_{j_{\sN}=0}^{\infty}\mQ^{\vec{\eta}_m}(\omega_{j_{\sN}}=u)\le m\frac{\beta}{d-1}G_{d-1}(0),\nn
    \enalign
where the last inequality holds by \refeq{Gbound1} with $i=0$.
\qed
\vskip0.5cm

{\noindent \em Proof of Proposition \ref{prp:pibound}}.
It follows from \refeq{piNxydef} and \refeq{Deltabound} that
\eqalign
\sum_{x,y}\sum_m|\pi_m^{\smallsup{N}}(x,y)|\le \Pi_{\sN}(\vec{f}_{\sN}),
\enalign
where the $\walkvec{i}{}$ for $i\ge0$ are excited random
walks and $\vec{f}_{\sN}$ is given by \refeq{fdefs}.

If $N=1$ then applying Lemma \ref{lem:decomposition} with $\kappa_1=1$,
$K_0=1$ and $K_1=\frac{\beta}{d}\mc{E}_0(d)$ (i.e. the right hand side
of the third bound of \refeq{fbounds}) we easily get the result.

For $N>1$, applying Lemma \ref{lem:decomposition} with
$\kappa_{\sN}=1$, $\kappa_i(j_i)=(j_i+1)$ for $i\ne N$, and
    \eqn{
    \lbeq{Ki-def}
    K_0=1, \quad K_1=\frac{\beta}{d}\mc{E}_1(d),\quad
    K_{\sN}=\frac{\beta}{d-1}G_{d-1}(0), \quad \text{and}\quad
    K_i=\beta a_d,~\text{for }2\leq i\leq N-1.
    }
(see the right hand sides of the remaining bounds of \refeq{fbounds}), we obtain
\refeq{pibound} for $N>1$.
\qed
\vskip0.5cm

Since the speed is known to exist \cite{BR07}, the following corollary is an easy consequence of
\cite[Propositions 3.1 and 6.1]{HH07}
together with Proposition \ref{prp:pibound}, and the fact that
$G_5^{*2}(0)<5^2/6$ \cite{HS92b}.
\begin{COR}[Formula for the speed of ERW]
\label{cor:ERWspeed}
For all $d\ge 6$ and $\beta \in [0,1]$,
    \eqalign
    \theta(\beta,d)=\lim_{n \ra \infty}\mE[\omega_{n+1}-\omega_n]
    =\frac{\beta e_1}{d}+\sum_{m=2}^{\infty}\sum_x x \pi_m(x).
    \enalign
\end{COR}
\ci{In fact, Corollary \ref{cor:ERWspeed} holds for all $d\ge 2$ since the law $\mu_n$ of the cookie environment as viewed by the random walker at time $n$ is known to converge (see e.g. \cite{BR07}).  Indeed
\eqalign
\mE[\omega_{n+1}-\omega_n]=\mE\big[\mE^{\vec{\omega}_n}[\omega_{n+1}-\omega_n]\big]=\mE\big[\frac{\beta e_1}{d}I_{\{\omega_n \notin \vec{\omega}_{n-1}\}}\big]=\frac{\beta e_1}{d}\big[1-\mP(\omega_n \in \vec{\omega}_{n-1})\big],
\enalign
where the right hand side converges as $n\ra \infty$ since $\mP(\omega_n \in \vec{\omega}_{n-1})$ is the $\mu_n$-measure of the event that the cookie at the origin is absent.  To prove monotonicity of the speed, it is therefore sufficient to prove that for each fixed $n$, $\mP_{\beta}(\omega_n \in \vec{\omega}_{n-1})$ is non-increasing in $\beta$.}

\section{The differentiation step}
\label{sec-derivpi}
To verify the exchange of limits in \refeq{deriv1},
it is sufficient to prove that $\sum_{x,y}(y-x)\pi_m^{\smallsup{N}}(x,y)$
is absolutely summable in $m$ and $N$ and that
$\sum_{N=1}^{\infty}\sum_{m=2}^\infty\sup_{\beta\in [0,1]}|\sum_{x,y}  (y-x)\varphi_m^{\smallsup{N}}(x,y)|<\infty$.
By Proposition \ref{prp:pibound} and the fact that $|y-x|=1$
for $x,y$ nearest neighbours, the first condition holds provided that
    \eqalign
    \boxed{\beta a_d<1.}\lbeq{C1}
    \enalign
In fact we will see later on that this inequality for $\beta=1$ is
sufficient to also establish the second condition. We now identify
$\varphi_m^{\smallsup{N}}(x,y)$.

For general $N$, we have (with $j_0=0$)
    \eqalign
    \pi_m^{\smallsup{N}}(x,y)=&\sum_{\vec{j}\in \mc{A}_{m,N}}\sum_{\walkvec{0}{1}}\sum_{\walkvec{1}{j_1+1}}\dots\sum_{\walkvec{N}{j_{N}+1} }I_{\{\omega^{(N)}_{j_{\sN}}=x, \omega^{(N)}_{j_{\sN}+1}=y\}}p^{\varnothing}(0,\walkcoor{0}{1})
    \prod_{n=1}^{N}\prod_{i_{n}=0}^{j_{n}-1}p^{\walkvec{n-1}{j_{n-1}+1}\circ \walkvec{n}{i_{n}}}\left(\walkcoor{n}{i_{n}},\walkcoor{n}{i_{n}+1}\right)
    \Delta^{\smallsup{n}}_{j_{n}+1}.
    \lbeq{pi-form}
    \enalign
Therefore,
    \eqalign
    \varphi_m^{\smallsup{N}}(x,y)=\varphi_m^{\smallsup{N,1}}(x,y)+\varphi_m^{\smallsup{N,2}}(x,y)+\varphi_m^{\smallsup{N,3}}(x,y),\lbeq{varphibreak}
    \enalign
where (by Leibniz' rule), $\varphi_m^{\smallsup{N,1}}(x,y)$, $\varphi_m^{\smallsup{N,2}}(x,y)$ and $\varphi_m^{\smallsup{N,3}}(x,y)$ arise from differentiating $p^{\varnothing}(0,\walkcoor{0}{1})$, $\prod_{n=1}^{N}\prod_{i_{n}=0}^{j_{n}-1}p^{\walkvec{n-1}{j_{n-1}+1}\circ \walkvec{n}{i_{n}}}\left(\walkcoor{n}{i_{n}},\walkcoor{n}{i_{n}+1}\right)$ and $\prod_{n=1}^{N}\prod_{i_{n}=0}^{j_{n}-1}\Delta^{\smallsup{n}}_{j_{n}+1}$, respectively.

Observe that if $\vec{\eta}_{m}=x_l$ then
    \eqalign
    \frac{\del}{\del\beta}{p_{\beta}^{\vec{\eta}_m}}(x_{l},x)=&\frac{e_1\cdot (x-x_l)I_{\{x_l \notin\vec{\eta}_{m-1}\}}}{2d}I_{\{|x-x_l|=1\}}=\frac{I_{\{x_l \notin\vec{\eta}_{m-1}\}}}{2d}\left(I_{\{x-x_l= e_1\}}-I_{\{x-x_l= -e_1\}}\right)\lbeq{pderiv},
    \enalign
and hence, using $I_{A}-I_{A\cap C}=I_{A \cap C^c}$ we have
    \eqalign
    \frac{\del}{\del\beta}\left({p_{\beta}^{\vec{\eta}_m}}(x_{l},x)-{p_{\beta}^{\vec{\omega}_n\circ\vec{\eta}_m}}(x_{l},x)\right)
    =\frac{1}{2d}I_{\{x_l \notin\vec{\eta}_{m-1},x_l \in \vec{\omega}_{n-1}\}}\left(I_{\{x-x_l= e_1\}}-I_{\{x-x_l= -e_1\}}\right).\lbeq{deltaderiv}
    \enalign
Clearly then
    \eqalign
    \left|\frac{\del}{\del\beta}\left({p_{\beta}^{\vec{\eta}_m}}(x_{l},x)-{p_{\beta}^{\vec{\omega}_n\circ\vec{\eta}_m}}(x_{l},x)\right)\right|
    \le \frac{1}{2d}I_{\{x_l \in \vec{\omega}_{n-1}\setminus \vec{\eta}_{m-1}\}}\left(I_{\{x-x_l= e_1\}}+I_{\{x-x_l= -e_1\}}\right).\lbeq{deltaderivrep}
    \enalign

Let $\rho^{\smallsup{N}}$ be obtained by replacing
$p^{\varnothing}(0,\walkcoor{0}{1})$ in \refeq{pi-form} with $(2d)^{-1}I_{\{\walkcoor{0}{1}= \pm e_1\}}$ (a bound on its derivative)
 and by bounding $\Delta^{\smallsup{n}}_{j_{n}+1}$ by $|\Delta^{\smallsup{n}}_{j_{n}+1}|$
for all $n=1, \ldots, N$.

For $k=1, \ldots, N$, let $\gamma^{\smallsup{N}}_k$
be obtained from \refeq{pi-form} by bounding
$\Delta^{\smallsup{n}}_{j_{n}+1}$ by $|\Delta^{\smallsup{n}}_{j_{n}+1}|$
for all $n=1, \ldots, N$ and by replacing $\prod_{i_{k}=0}^{j_{k}-1}p^{\walkvec{n-1}{j_{n-1}+1}\circ \walkvec{k}{i_{k}}}\left(\walkcoor{k}{i_{k}},\walkcoor{k}{i_{k}+1}\right)$ with the following bound on its derivative
\eq
\sum_{l=0}^{j_{k}-1}\frac{I_{\{\walkcoor{k}{l_{k}+1}-\walkcoor{k}{l_{k}}= \pm e_1\}}}{2d}\prodtwo{i_{k}=0}{i_{k}\ne l}^{j_{k}-1}p^{\walkvec{n-1}{j_{n-1}+1}\circ \walkvec{k}{i_{k}}}\left(\walkcoor{k}{i_{k}},\walkcoor{k}{i_{k}+1}\right).
\en
Similarly, let $\gamma^{\smallsup{N}}_k$
be obtained by replacing $\Delta^{\smallsup{k}}_{j_{k}+1}$ in
\refeq{pi-form} by
$(2d)^{-1}I_{\{\walkcoor{k}{j_{k}} \in \walkvec{k-1}{j_{k-1}+1}\setminus\walkvec{k}{j_{k}-1}\}}
I_{\{\walkcoor{k}{j_{k}+1}-\walkcoor{k}{j_{k}}= \pm e_1\}}$ (a bound on its derivative) and by bounding
$\Delta^{\smallsup{n}}_{j_{n}+1}$ for $n\neq k$ by $|\Delta^{\smallsup{n}}_{j_{n}+1}|$.

Letting $\gamma^{\smallsup{N}}=\sum_{k=1}^{N}\gamma^{\smallsup{N}}_k$ and
$\chi^{\smallsup{N}}=\sum_{k=1}^{N}\chi^{\smallsup{N}}_k$, we obtain that
    \eqalign
    \sum_{m}\sum_{x,y}|\varphi_m^{\smallsup{N,1}}(x,y)|\le \rho^{\smallsup{N}}, \quad \sum_{m}\sum_{x,y}|\varphi_m^{\smallsup{N,2}}(x,y)|\le
    \gamma^{\smallsup{N}}, \quad \text{and  }\sum_{m}\sum_{x,y}|\varphi_m^{\smallsup{N,3}}(x,y)|\le
    \chi^{\smallsup{N}}.\lbeq{3terms}
    \enalign

\begin{LEM}[Bounds on $\rho^{\smallsup{N}}$]
\label{lem:rhobound}
We have $\rho^{\smallsup{1}}\le d^{-2}\beta\mc{E}_0(d)$, and, for $N\ge 2$,
    \eqalign
    \rho^{\smallsup{N}}\le \beta^N
    \frac{G_{d-1}(0)\mc{E}_1(d)}{d^2(d-1)}
    a_d^{N-2}.
    \lbeq{rhobound}
    \enalign
\end{LEM}

\proof For $N\ge 1$,
    \eqalign
    \rho^{\smallsup{N}}&=\frac{1}{d}\sum_{\vec{j}\in \mc{A}_{N}}\sum_{\walkvec{0}{j_0+1}}\mE^{\dagger\varnothing}\Bigg[\mE^{\walkvec{0}{j_0+1}}\Big[\sum_{\walkcoor{1}{j_1+1}}|\Delta^{\smallsup{1}}_{j_{1}+1}|\dots \nn\\
    &\qquad \mE^{\walkvec{N-2}{j_{N-2}+1}}\big[\sum_{\walkcoor{N-1}{j_{N-1}+1}}|\Delta^{\smallsup{N-1}}_{j_{N-1}+1}|
    \mE^{\walkvec{N-1}{j_{N-1}+1}}\big[\sum_{\walkcoor{N}{j_{\sN}+1}}|\Delta^{\smallsup{N}}_{j_{N}+1}|\big]\big]\dots\Big]
    \Bigg]\le \Pi_{\sN}(\vec{g}_{\sN})\lbeq{summrho},
    \enalign
where $g_{0,j_0}=d^{-1}f_{0,j_0}$, $g_{i,j_i}=f_{i,j_i}$ for $i\ge 1$,
and the $\walkvec{i}{}$ for $i\ge1$ are excited random walks, while
$\walkvec{0}{}$ is a 1-step simple random walk in the {\em first coordinate only}.
This latter difference is indicated by the dagger in the notation
$\mE^{\dagger\varnothing}$.  Since we have already established the
relevant bounds on the $f_{i,j_i}$, to complete the proof of
Lemma \ref{lem:rhobound} by applying Lemma \ref{lem:decomposition},
it is enough to establish that
    \eqalign
    \sum_{j_0=0}^{\infty}(j_0+1)\mE^{\dagger\vec{\eta}}_0[g_{0,{j}_0}]\le \frac{1}{d},\nn
    \enalign
which is trivial.
\qed

\begin{LEM}[Bounds on $\chi^{\smallsup{N}}$]
\label{lem:chibound}
We have $\chi^{\smallsup{1}}\leq d^{-1}\mc{E}_0(d)$, and, for $N\ge 2$,
    \eqalign
    \chi^{\smallsup{N}}\le
    N\beta^{N-1}\frac{G_{d-1}(0)\mc{E}_1(d)}{d(d-1)}
    a_d^{N-2}.
    \lbeq{chibound}
    \enalign
\end{LEM}
\proof We rewrite
    \eqalign
    \chi^{\smallsup{N}}_k&=\sum_{\vec{j}\in \mc{A}_{N}}\sum_{\walkvec{0}{j_0+1}}\mE^{\varnothing}\Bigg[\mE^{\walkvec{0}{j_0+1}}\Big[\sum_{\walkcoor{1}{j_1+1}}|\Delta^{\smallsup{1}}_{j_{1}+1}|\dots \mE^{\walkvec{k-1}{j_{k-1}+1}}\Big[\sum_{\walkcoor{k}{j_{k}+1}}\frac{I_{\{
    \walkcoor{k}{j_{k}} \in \walkvec{k-1}{j_{k-1}+1}\setminus\walkvec{k}{j_{k}-1}\}}}{2d}I_{\{\walkcoor{k}{j_{k}+1}-\walkcoor{k}{j_{k}}= \pm e_1\}}\dots\nn\\
    &\qquad \mE^{\walkvec{N-2}{j_{N-2}+1}}\big[\sum_{\walkcoor{N-1}{j_{N-1}+1}}
    |\Delta^{\smallsup{N-1}}_{j_{N-1}+1}|\mE^{\walkvec{N-1}{j_{N-1}+1}}\big[\sum_{\walkcoor{N}{j_{\sN}+1}}
    |\Delta^{\smallsup{N}}_{j_{N}+1}|\big]\big]\dots\Big]\Bigg]
    \le\Pi_{\sN}(\vec{\phi}^{\sss(k)}_{\sN}),
    \enalign
where $\phi^{\sss(k)}_{i,j_i}=f_{i,j_i}$ for $i\ne k$, and $\phi^{\sss(k)}_{k,j_k}=\beta^{-1}f_{k,j_k}$, and the $\walkvec{i}{}$ for $i\ge0$ are excited random walks.  The resulting bound on $\chi^{\smallsup{N}}$, (which is simply $\beta^{-1}N$ times \refeq{pibound}) is then easily obtained by applying Lemma \ref{lem:decomposition} to each of the $\chi^{\smallsup{N}}_k$ and summing over $k$.
\qed

Before proceeding to the bound on $\gamma^{(N)}$, we first need
a new lemma similar to Lemma \ref{lem:togreens}.

\begin{LEM}[Ingredients for bounds on the {\em derivative} of the speed of ERW]
\label{lem:togreens2}
Let $\mathbb{Q}^{\leftrightarrow_l,\vec{\eta}_m}$ denote the law of a
self-interacting random walk with history $\vec{\eta}_m$, where the walk
is an excited random walk, except for the $l^{\rm th}$-step, which is a
simple random walk step in the first coordinate only.  Then, for $i\geq 0$,
    \eqalign
    \sum_{j=1}^{\infty}\frac{(j+i)!}{j!}\sum_{l=1}^{j}\mQ^{\leftrightarrow_l,\vec{\eta}_m}(\omega_{j}=u)\le &(i+1)!\left(\frac{d}{d-1}\right)^{i+2}G_{d-1}^{*(i+2)}(0).\lbeq{greensb2}
    \enalign
\end{LEM}

\proof Since one of the $j$ steps is a simple random walk step in the
first coordinate, the number of steps in the other coordinates has a
Binomial$(j-1, \frac{d-1}{d})$ distribution.
Thus,
    \eqalign
    \sum_{j=1}^{\infty}\frac{(j+i)!}{j!}\sum_{l=1}^{j}\mQ^{\leftrightarrow_l,\vec{\eta}_m}(\omega_{j}=u)\le
    &\sup_{v\in \Z^{d-1}} \sum_{j=1}^{\infty}\sum_{l=1}^{j}\sum_{k=0}^{j-1}\frac{(j+i)!}{j!}\mc{P}(\mc{N}_{j-1}=k)\mP_{d-1}(\omega_k=v)\nn\\
    =&\sup_{v\in \Z^{d-1}} \sum_{k=0}^{\infty}\mP_{d-1}(\omega_k=v)\sum_{j=k+1}^{\infty}j\frac{(j+i)!}{j!}\mc{P}(\mc{N}_{j-1}=k)\nn\\
    =&\sup_{v\in \Z^{d-1}} \sum_{k=0}^{\infty}\mP_{d-1}(\omega_k=v)\sum_{r=k}^{\infty}\frac{(r+i+1)!}{r!}
    \mc{P}(\mc{N}_{r}=k).
    \enalign
Now proceed as in the proof of Lemma \ref{lem:togreens} to obtain the result.
\qed


Define
    \eqn{
    \epsilon(d)=
    \frac{2d}{(d-1)^4}G_{d-1}(0)G_{d-1}^{*3}(0)+\frac{\mc{E}_1(d)}{d(d-1)^2}G_{d-1}^{*2}(0),
    \lbeq{epsilonddef}
    }

\begin{LEM}[Bounds on $\gamma^{\smallsup{N}}$]
\label{lem:gammabound}
We have $\gamma^{\smallsup{1}}\le \beta(d-1)^{-2}G_{d-1}^{*2}(0)$,
$\gamma^{\smallsup{2}}\le\beta^2\epsilon(d)$ and, for all $N\ge 3$,
    \eqn{
    \gamma^{\smallsup{N}}\le
    \epsilon(d)\beta^2(\beta a_d)^{N-2}+
    (N-2)\frac{2\beta^3\mc{E}_1(d)}{(d-1)^4}G_{d-1}(0)G_{d-1}^{*3}(0)
    (\beta a_d)^{N-3}.
    }
\end{LEM}
\proof We rewrite
    \eqalign
    \gamma^{\smallsup{N}}_k&=\sum_{\vec{j}\in \mc{A}_{N}}\mE^{\varnothing}\Bigg[\mE^{\walkvec{0}{1}}\Bigg[\sum_{\walkcoor{1}{j_1+1}}|\Delta^{\smallsup{1}}_{j_{1}+1}|
    \dots \mE^{\walkvec{k-2}{j_{k-2}+1}}\Big[\sum_{\walkcoor{k-1}{j_{k-1}+1}}|\Delta^{\smallsup{k-1}}_{j_{k-1}+1}|\nn\\
    &\qquad \sum_{l=0}^{j_{k}-1}\mE^{\leftrightarrow_l \walkvec{k-1}{j_{k-1}+1}}\Big[\sum_{\walkcoor{k}{j_{k}+1}}|\Delta^{\smallsup{k}}_{j_{k}+1}|
    \mE^{\walkvec{k}{j_{k}+1}}\big[\sum_{\walkcoor{k+1}{j_{k+1}+1}}|\Delta^{\smallsup{k+1}}_{j_{k+1}+1}|\dots \mE^{\walkvec{N-1}{j_{N-1}+1}}\big[\sum_{\walkcoor{N}{j_{\sN}+1}}|\Delta^{\smallsup{N}}_{j_{N}+1}|\big]\dots\big]\dots\Big]\Big]\Bigg]\Bigg]\nn\\
    &\leq \Pi_{\sN}^{\sss(k)}(\vec{h}_{\sN})\lbeq{pigamma},
    \enalign
where $h_{i,j_i}=f_{i,j_i}$ for $i\ne k$ and $h_{k,j_k}=d^{-1}f_{k,j_k}$,
and the $\walkvec{i}{}$ for $i\ne k$ are excited random walks, while
$\walkvec{k}{}$ is an excited random walk except that its $(l+1)$st
step is a simple random walk step in the first coordinate.  This is
indicated by the left-right arrow with subscript $l$ in the
notation $\mE^{\leftrightarrow_l \walkvec{k-1}{j_{k-1}+1}}$.

When $N=1$, then also $k=1$ and we use \refeq{greensb2} with $i=0$, and
Lemma \ref{lem:decomposition} to get the required bound.

When $N>1$ and $k=1$ we use the same bounds as in the proof
of Proposition \ref{prp:pibound} except that we use
\refeq{greensb2} with $i=1$ on the term $k=1$.  This gives us a
bound on $\gamma^{\smallsup{N}}_1$ (when $N>1$) of
    \eqalign
    \frac{2\beta}{d^2}\left(\frac{d}{d-1}\right)^3
    G_{d-1}^{*3}(0)
    \frac{\beta}{d-1}G_{d-1}(0)
    \prod_{i=2}^{N-1}\beta a_d.\lbeq{kis1}
    \enalign
When $N>1$ and $k=N$, we use the same bounds as in the
proof of Proposition \ref{prp:pibound} except that we use
\refeq{greensb2} with $i=0$ on the term $k=N$.  This gives
us a bound on $\gamma^{\smallsup{N}}_{\sN}$ (when $N>1$) of
    \eqalign
    \frac{\beta}{(d-1)^2}G_{d-1}^{*2}(0)\frac{\beta}{d}\mc{E}_1(d)
    \prod_{i=2}^{N-1}\beta a_d.\lbeq{kisN}
    \enalign
Similarly when $N>1$ and $1\ne k\ne N$ (so $N>2$) we use \refeq{greensb2} on term $k$ to get a bound on $\gamma^{\smallsup{N}}_k$ of
    \eqalign
    \frac{\beta}{d-1}G_{d-1}(0)
    \frac{\beta}{d}\mc{E}_1(d)
    \frac{2\beta}{d^2}\left(\frac{d}{d-1}\right)^3G_{d-1}^{*3}(0)
    \prodstack{i=2}{i\ne k}^{N-1}\beta a_d.\lbeq{kisN2}
    \enalign
Simplifying these expressions and summing over $k$ completes the proof of the lemma. \qed

\begin{COR}[Summary of bounds]
\label{cor:bounds}
For all $\beta \in [0,1]$, and $d$ such that $a_d<1$
    \eqalign
    d\sum_{N=1}^{\infty}\rho^{(N)}\le
    &\frac{\mc{E}_0(d)}{d}+\frac{G_{d-1}(0)\mc{E}_1(d)}{d(d-1)(1-a_d)}\lbeq{RHO}\\
    d\sum_{N=1}^{\infty}\chi^{(N)}\le &\mc{E}_0(d)+\frac{G_{d-1}(0)\mc{E}_1(d)(2-a_d)}{(d-1)(1-a_d)^2}\lbeq{CHI}\\
    d\sum_{N=1}^{\infty}\gamma^{(N)}\le & \frac{dG_{d-1}^{*2}(0)}{(d-1)^2}+\frac{\epsilon(d)d}{1-a_d}+
    \frac{2d\mc{E}_1(d)G_{d-1}(0)G_{d-1}^{*3}(0)}{(d-1)^4(1-a_d)^{2}}.\lbeq{GAMMA}
    \enalign
\end{COR}
\proof Firstly note that the condition on $a_d$ ensures that $\rho^{\sss (N)}$, $\chi^{\sss (N)}$
and $\gamma^{\sss(N)}$ are all summable over $N$, and in all cases the supremum
over $\beta$ occurs at $\beta=1$ (see Lemmas \ref{lem:rhobound},
\ref{lem:chibound} and \ref{lem:gammabound}).  The results are then
easily obtained by summing each of the bounds in Lemmas
\ref{lem:rhobound}, \ref{lem:chibound} and \ref{lem:gammabound}
over $N$.
\qed

\section{Proof of Theorem \ref{thm:main}}
\label{sec-pfmainthm}
For $d$ such that $a_d<1$, the bounds of Corollary \ref{cor:bounds}
hold.  From \refeq{3terms} we have the required absolute summability
conditions in the discussion after \refeq{deriv1}, and in particular
\refeq{deriv1} holds for all $\beta$.  To complete the proof of the
theorem, it remains to show that the right hand side of \refeq{needed1}
is no more than $d^{-1}$.  By \refeq{3terms} and Corollary \ref{cor:bounds},
we have bounded $d$ times the right hand side of \refeq{needed1} by the sum
of the right hand sides of the bounds in Corollary \ref{cor:bounds}.
Since these terms all involve simple random walk Green's functions
quantities, we will need to use estimates of these quantities.

In order to bound $\mc{E}_i(d)$, we shall first prove that,
for all $i\geq 0$,
    \eqn{
    \lbeq{Eirewrfin}
    \mc{E}_i(d)=
    \big(\frac{d}{d-1}\big)^{i+1}G_{d-1}^{*(i+1)}(0)-1.
    }
In order to prove \refeq{Eirewrfin}, we first make use of
\cite[Lemma B.3]{HS92b}, which states that $G_d^{*n}(x)$ is
non-increasing in $|x_i|$ for every $i=1, \ldots, d$, so that
the supremum in \refeq{Eidef} can be restricted to $v=0$ and
$v=e$ for any neighbour $e$ of the origin. In order to bound
$G_d^{*n}(e)$, we make use of the fact that for any
function $x\mapsto f(x)$ for which $f(e)$ is constant for all
$e\in \Z^d$ with $|e|=1$, we have $f(e)=(D_d*f)(0)$, so that
    \eqalign
    \mc{E}_i(d)=&
    \max\Big\{\big(\frac{d}{d-1}\big)^{i+1}G_{d-1}^{*(i+1)}(0)-1, \big(\frac{d}{d-1}\big)^{i+1}(D_{d-1}*G_{d-1}^{*(i+1)})(0)\Big\}.
    \lbeq{Eirewr}
    \enalign
Finally, note that since $G_{d}(x)=\delta_{0,x}+(D_d* G_d)(x)$, we have
that $G_{d}^{*(i+1)}(0)=G_{d}^{*i}(0)+(D_d* G_d^{*(i+1)})(0)$.
Therefore, since $G_{d}^{*i}(0)\geq 1$ for all $i\geq 0$,
    \eqalign
    \big(\frac{d}{d-1}\big)^{i+1}G_{d-1}^{*(i+1)}(0)-1
    &=\big(\frac{d}{d-1}\big)^{i+1}(D_{d-1}*G_{d-1}^{*(i+1)})(0)
    +\big(\frac{d}{d-1}\big)^{i+1}G_{d}^{*i}(0)-1\nn\\
    &>\big(\frac{d}{d-1}\big)^{i+1}(D_{d-1}*G_{d-1}^{*(i+1)})(0),
    \enalign
which proves \refeq{Eirewrfin}.

By \cite[Lemma C.1]{HS92b}, $d\mapsto G_{d}^{*n}(0)$
is monotone decreasing in $d$ for each $n\geq 1$, so that it
suffices to show that the sum of terms on the
right hand sides of \refeq{RHO}, \refeq{CHI} and \refeq{GAMMA}
is bounded by $1$ for $d=\dimres$.
For this we use the following rigorous Green's
functions estimates \cite{Hpc, HS92b} \ci{for $d=8$:
    \ci{\eqalign
    &G_d(0)\le 1.07865, \quad G^{*2}_d(0)\le 1.2891, \quad G^{*3}_d(0)\le 1.8316.
    \enalign
    }
Putting in these values for $d-1=8$ we get that the sum of the
right hand sides of the bounds in Corollary \ref{cor:bounds}
is at most $0.97$, whence the result follows for $d\geq \dimres$.
%
}
\qed

To prove monotonicity for $\beta \in [0,\beta_0]$ for
some $\beta_0(d)$ for each $d\ge \dimlow$, it is sufficient to
prove that $\chi^{\sss (1)}<d^{-1}$ when $d\ge \dimlow$
(and that the other terms are bounded), since
this is the only term that does not contain the small
factor $\beta$.  Since $\chi^{\sss (1)}\le d^{-1}\mc{E}_0(d)$,
it is enough to show that $\mc{E}_0(d)<1$ for $d=\dimlow$,
since the right hand sides of \refeq{RHO}, \refeq{CHI}
and \refeq{GAMMA} are bounded for $d\geq \dimlow$.
From \cite{HS92b} we have $\frac{6}{5}G_5(0)-1<\frac{6}{5}(1.157)-1<1$, and since $\mc{E}_0(d)$ is decreasing in $d$, this completes the result.




\paragraph{Acknowledgements.}
The work of RvdH and MH was supported in part by Netherlands Organisation for
Scientific Research (NWO).  The work of MH was also supported by a FRDF
grant from the University of Auckland.  The authors would like to thank Itai Benjamini for suggesting this problem to us and 
Takashi Hara for providing the SRW Green's functions upper bounds.

\bibliographystyle{plain}



\end{document}

%% file: def99c.tex
\def\UseSection{
        \numberwithin{equation}{section}
	\theoremstyle{plain}
        \newtheorem{theorem}    {Theorem}[section]
        \DefineTheorems 
}

\def\DefineTheorems{
	
	\newtheorem{lemma}      [theorem] {Lemma}
	
	\newtheorem{prop}       [theorem] {Proposition}
	
	\newtheorem{cor}        [theorem] {Corollary}

	\theoremstyle{definition}
	\newtheorem{defn}       [theorem] {Definition}

	\theoremstyle{definition}

}

\newcommand{\bt}   {\begin{theorem}}
\newcommand{\et}   {\end  {theorem}}
\newcommand{\bl}   {\begin{lemma}}
\newcommand{\el}   {\end  {lemma}}
\newcommand{\bp}   {\begin{prop}}
\newcommand{\ep}   {\end  {prop}}
\newcommand{\bc}   {\begin{cor}}
\newcommand{\ec}   {\end  {cor}}
\newcommand{\bd}   {\begin{defn}}
\newcommand{\ed}   {\end  {defn}}

\newcommand{\ba}   {\begin{array}}
\newcommand{\ea}   {\end  {array}}
\newcommand{\be}   {\begin{enumerate}}
\newcommand{\ee}   {\end  {enumerate}}
\newcommand{\bi}   {\begin{itemize}}
\newcommand{\ei}   {\end  {itemize}}

\def\eq#1\en{\begin{equation}#1\end{equation}}  
\def\eqsplit#1\ensplit{
	\begin{equation}\begin{split}#1\end{split}\end{equation}
	}
\def\eqalign#1\enalign{
	\begin{align}#1\end{align}
	}
\def\eqmul#1\enmul{
	\begin{multline}#1\end{multline}
	}
\newcommand{\eqarrstar} {\begin{eqnarray*}} 
\newcommand{\enarrstar} {\end{eqnarray*}} 
\newcommand{\eqarray}   {\begin{eqnarray}} 
\newcommand{\enarray}   {\end{eqnarray}}

\newcommand{\lbeq}[1]  {\label{e:#1}}
\newcommand{\refeq}[1] {\eqref{e:#1}}    

\makeatletter
\newcommand{\labelcounter}[2]{{%
	\stepcounter{#1}
	\protected@write\@auxout{}%
	{\string\newlabel{#2}{{\csname the#1\endcsname}{\thepage}}}%
	{\ref{#2}}
	}}
\makeatother
\newcommand{\sss}   { \scriptscriptstyle } 

\newcommand{\Ebold} {{\mathbb E}}

\newcommand{\Qbold} {{\mathbb Q}}

\newcommand{\Zbold} {{\mathbb Z}}

\newcommand{\Zd}    {{ {\Zbold}^d }}

\newcommand{\spose}[1] {{\hbox to 0pt{#1\hss}} }
\newcommand{\ltapprox} {\mathrel{\spose{\lower 3pt\hbox{$\mathchar"218$}}
 \raise 2.0pt\hbox{$\mathchar"13C$}}}
\newcommand{\gtapprox} {\mathrel{\spose{\lower 3pt\hbox{$\mathchar"218$}}
 \raise 2.0pt\hbox{$\mathchar"13E$}}}

\newcommand{\prodtwo}[2]{
	\prod_{ \mbox{ \scriptsize 
		$\begin{array}{c} 
		{#1} \\ 
		{#2}  
		\end{array} $ } 
		} 
	} 




%% file: mono_fin.bbl
\begin{thebibliography}{1}

\bibitem{AR05}
T.~Antal and S.~Redner.
\newblock The excited random walk in one dimension.
\newblock {\em J. of Physics A.}, {\bf 38}:2555--2577, (2005).

\bibitem{BW03}
I.~Benjamini and D.~B. Wilson.
\newblock Excited random walk.
\newblock {\em Electron. Comm. Probab.},{\bf 8}:86--92 (electronic), (2003).

\bibitem{BR07}
J.~B{\'e}rard and A.~Ram{\'i}rez.
\newblock {Central limit theorem for excited random walk in dimension $d\ge
  2$}.
\newblock {\em Electronic Communications in Probability}, {\bf 12}: 303--314, (2007).

\bibitem{BolSznZei03}
E.~Bolthausen, A.-S. Sznitman, and O.~Zeitouni.
\newblock Cut points and diffusive random walks in random environment.
\newblock {\em Ann. Inst. H. Poincar\'e Probab. Statist.}, {\bf
  39}(3):527--555, (2003).

\bibitem{Dav99}
B.~Davis.
\newblock Brownian motion and random walk perturbed at extrema.
\newblock {\em Probab. Theory Relat. Fields.}, {\bf 113}:501--518, (1999).

\bibitem{Hpc}
T.~Hara.
\newblock Private communication, (2007).


\bibitem{HS92b}
T.~Hara and G.~Slade.
\newblock The lace expansion for self-avoiding walk in five or more dimensions.
\newblock {\em Reviews in Math.\ Phys.}, {\bf 4}:235--327, (1992).

\bibitem{HH07}
R.~van~der Hofstad and M.~Holmes.
\newblock An expansion for self-interacting random walks.
\newblock arXiv:0706.0614v3 [math.PR], (2007).

\bibitem{Kozm03}
G.~Kozma.
\newblock Excited random walk in three dimensions has positive speed.
\newblock (2003).

\bibitem{Kozm05}
G.~Kozma.
\newblock Excited random walk in two dimensions has linear speed.
\newblock arXiv:math/0512535v1 [math.PR], (2005).

\bibitem{SZ99}
A.-S.~Sznitman and M.~Zerner
\newblock A law of large numbers for random walks in random environment.
\newblock {\em Ann. Probab.}, {\bf 27}: 1851 - 1869, (1999).

\bibitem{Zern05}
M.~Zerner.
\newblock Multi-excited random walks on integers.
\newblock {\em Probab. Theory Relat. Fields.}, {\bf 133}:98--122, (2005).

\end{thebibliography}
